# Fractional Quantization of Singular Lagrangian systems with Second-Order Derivatives Using WKB Approximation

Eyad Hasan Hasan

Tafila Technical University, Faculty of Science, Applied Physics Department, P. O. Box: 179, Tafila 66110, Jordan

## Abstract

In this paper, the theory of the fractional singular Lagrangian systems is investigated with second-order derivatives. The fractional quantization for these systems is examined using the WKB approximation. The Hamilton–Jacobi treatment can be applied for these systems. The equations of motion are obtained. The Hamilton–Jacobi partial differential equations are solved to determine the action function. By finding the action function, the wave function for these systems is constructed. We achieved that the quantum results are agreement with the classical results. Besides, to demonstrate the theory; two mathematical examples are examined.

CONTENTS

**Email**: iyad973@yahoo.com, dr_eyad2004@ttu.edu.jo

## 1. Introduction

The quantization of singular Lagrangian systems have been treated with more interests by Dirac's work for quantizing the gravitational field [1, 2]. Following Dirac's work, researchers developed a new formalism for investigating singular Lagrangian system which is called the canonical method was developed by Guler [3], then this formalism has been used to quantize these systems using the WKB approximation and path integral approach [4-8]. A general theory has been investigated for quantizing higher-order singular Lagrangian systems using WKB approximation by Hasan et al [6-8]. In this theory, researchers have achieved that the quantum results are found to be in exact agreement with the classical results.

In this paper, we would like to develop this theory to be applicable for these systems within fractional derivatives. The quantization of fractional singular Lagranians systems have been studied with increasing interest for quantizing physical systems [9, 10]. Recently, fractional Lagrangians systems with second-order derivatives have been treated with more interests and importance [11, 12]. Researchers have investigated the Hamilton-Jacobi formalism for these systems within fractional derivatives, the Euler-Lagrange equations and Hamilton's equations are analyzed [11, 12]. More recently, researchers have constructed a formalism using the canonical method [3] for quantization singular systems using the WKB approximation and path integral approach for first-order derivatives [9, 10]. In this formalism, the equations of motion are written in fractional form as total differential equations, also, the set of fractional Hamilton–Jacobi partial differential equations is constructed. In addition, the solution of this set HJPDEs in fractional form and the fractional Hamilton-Jacobi function are obtained.

In this paper, we would like to extend the work for fractional singular Lagrangians systems with second-order derivatives and quantize these systems using WKB approximation. We constructed a formalism for investigating fractional singular Lagrangian systems and Hamilton-Jacobi formalism for second-order derivatives. Besides, we constructed the fractional Hamilton-Jacobi function for these systems. By finding this function enable us to build the appropriate wave function, then, the quantization for these systems using WKB approximation can be carried out.

Now, the most important definitions of fractional derivatives are the left Riemann–Liouville and right Riemann–Liouville are defined respectively as [13].

$$ {}_aD_t^{\alpha} f(t) = \frac{1}{\Gamma(n-\alpha)} \left( \frac{d}{dt} \right)^n \int_a^t (t-\tau)^{n-\alpha-1} f(\tau) d\tau \,. \tag{1} $$

$$ {}_tD_b^{\alpha} f(t) = \frac{1}{\Gamma(n-\alpha)} \left( -\frac{d}{dt} \right)^n \int_t^b (\tau-t)^{n-\alpha-1} f(\tau) d\tau \,. \tag{2} $$

where $n \in N$, $n-1 \le \alpha < n$ and $\Gamma$ is the Euler gamma function. Here $\alpha$ is an integer, and these derivatives have properties as follows:

$$ {}_aD_t^{\alpha} f(t) = \left( \frac{d}{dt} \right)^{\alpha} f(t) \,, \tag{3} $$

$$ {}_tD_b^{\alpha} f(t) = \left( -\frac{d}{dt} \right)^{\alpha} f(t) \,, \tag{4} $$

In this work, we aim to construct the formalism for quantizing singular Lagrangians systems with second-order derivatives within framework of fractional derivatives. This paper is organized as follows: In section 2, we constructed the fractional singular Lagrangian and Hamilton-Jacobi formalism for second-order derivatives. In section 3, the fractional Hamilton-Jacobi function and quantization using WKB approximation are discussed. In section 4, Two illustrative examples are examined. The work closes with some concluding remarks in section 5.

## 2. Fractional Singular Lagrangian and Hamilton-Jacobi Formalism for Second-Order Derivatives

In this section, we will formulate the second-order singular Lagrangian within framework of fractional calculus. The starting point is a Lagrangian depending on the fractional derivatives is given by

$$ L = L(D^{\alpha-1} q_i, D^{\alpha} q_i, D^{2\alpha} q_i, t) \,. \tag{5} $$

Thus , the fractional of the Hessian matrix is defined as

$$W_{ij} = \frac{\partial^2 L}{\partial D^{2\alpha} q_i \partial D^{2\alpha} q_j} \qquad i,j = 1,2,\ldots,N \tag{6}$$

If it's rank is $N$, the fractional Lagrangian is called regular otherwise the Lagrangian is singular $N-R, \quad R<N$. Now, we can define the generalized momenta $\pi_i$ corresponding to the generalized coordinates $D^{\alpha} q_i$ as:

$$\pi_a = \frac{\partial L}{\partial D^{2\alpha} q_a}, \tag{7}$$

$$\pi_\mu = \frac{\partial L}{\partial D^{2\alpha} q_\mu}. \tag{8}$$

Since the rank of the Hessian matrix is N-R, one can solve Eq. (7) to obtain N-R accelerations $D^{2\alpha} q_a$ in terms of $D^{\alpha-1} q_i, D^{\alpha} q_i, \pi_a$ and $D^{2\alpha} q_\mu$ as follows:

$$D^{2\alpha} q_a = w_a (D^{\alpha-1} q_i, D^{\alpha} q_i, \pi_a, D^{2\alpha} q_\mu). \tag{9}$$

Substituting (9) in (8), we can obtain

$$\pi_\mu = -H_\mu^\pi (D^{\alpha-1} q_i, D^{\alpha} q_i, p_a, \pi_a) \tag{10}$$

A similar expression for the momenta $p_\mu$ can be obtained as:

$$p_\mu = -H_\mu^p (D^{\alpha-1} q_i, D^{\alpha} q_i, p_a, \pi_a) \tag{11}$$

and the generalized momenta $p_i$ corresponding to the generalized coordinate $D^{\alpha-1} q_i$ can be written as:

$$p_a = \frac{\partial L}{\partial D^{\alpha} q_a} - \frac{d}{dt}\left( \frac{\partial L}{\partial D^{2\alpha} q_a} \right); \tag{12a}$$

$$p_\mu = \frac{\partial L}{\partial D^{\alpha} q_\mu} - \frac{d}{dt}\left( \frac{\partial L}{\partial D^{2\alpha} q_\mu} \right). \tag{12b}$$

$a = 1.2,\ldots,N-R$; $\mu = 1,\ldots,R$.

Equations (10) and (11) can be written as

$$H'^{p}_{\mu}(D^{\alpha-1}q_i, D^{\alpha}q_i, p_i, \pi_i) = p_\mu + H^{p}_{\mu} = 0 \quad ; \qquad (13a)$$

$$H'^{\pi}_{\mu}(D^{\alpha-1}q_i, D^{\alpha}q_i, p_i, \pi_i) = \pi_\mu + H^{\pi}_{\mu} = 0 \qquad (13b)$$

Thus, equations (13) represent primary constraints [1, 2].

Now, we will define the Poisson bracket

$$\{A, B\} \equiv \frac{\partial A}{\partial D^{\alpha-1}q_i}\frac{\partial B}{\partial p_i} - \frac{\partial A}{\partial p_i}\frac{\partial B}{\partial D^{\alpha-1}q_i} + \frac{\partial A}{\partial D^{\alpha}q_i}\frac{\partial B}{\partial \pi_i} - \frac{\partial A}{\partial \pi_i}\frac{\partial B}{\partial D^{\alpha}q_i} .$$

Where A and B are functions in the phase space described in term of the canonical variables $D^{\alpha-1}q_i$, $D^{\alpha}q_i$, $p_i$ and $\pi_i$. Here, the generalized momenta $p_i$ and $\pi_i$ are conjugated to the generalized coordinates $D^{\alpha-1}q_i$ and $D^{\alpha}q_i$ respectively. Thus, the fundamental Poisson brackets are

$$\{D^{\alpha-1}q_i, p_j\} \equiv \delta_{ij} \text{ and } \{D^{\alpha}q_i, \pi_j\} \equiv \delta_{ij}, \qquad \text{where } i, j = 1, \ldots, N$$

The fractional Hamiltonian $H_\circ$ can be defined as

$$H_\circ = -L(D^{\alpha-1}q_i, D^{2\alpha}q_v, D^{\alpha}q_v, W_a) + p_a D^{\alpha}q_a + \pi_a D^{2\alpha}q_a - D^{\alpha}q_\mu H^{p}_{\mu} - D^{2\alpha}q_\mu H^{\pi}_{\mu} . \quad (14)$$

$$\mu = 1, \ldots\ldots, R\ ; \qquad a = R+1, \ldots, N\ .;$$

Because of a natural of singular Lagrangian, the generalized momenta $p_\mu$ and $\pi_\mu$ are not independent of $p_a$ and $\pi_a$. Thus, we can write the set of Hamilton-Jacobi partial differential equations as [14]

$$H'_{\circ}(D^{\alpha-1}q_i, D^{\alpha}q_i, \frac{\partial S}{\partial D^{\alpha-1}q_a}, \frac{\partial S}{\partial D^{\alpha-1}q_v}, \frac{\partial S}{\partial D^{\alpha}q_a}, \frac{\partial S}{\partial D^{\alpha}q_v}\Bigg) = p_{\circ} + H_{\circ} = 0; \quad (15a)$$

$$H'^{p}_{\mu}(D^{\alpha-1}q_i, D^{\alpha}q_i, \frac{\partial S}{\partial D^{\alpha-1}q_a}, \frac{\partial S}{\partial D^{\alpha-1}q_v}, \frac{\partial S}{\partial D^{\alpha}q_a}, \frac{\partial S}{\partial D^{\alpha}q_v}\Bigg) = p_{\mu} + H^{p}_{\mu} = 0 \ ; \quad (15b)$$

$$H'^{p}_{\mu}(D^{\alpha-1}q_i, D^{\alpha}q_i, \frac{\partial S}{\partial D^{\alpha-1}q_a}, \frac{\partial S}{\partial D^{\alpha-1}q_v}, \frac{\partial S}{\partial D^{\alpha}q_a}, \frac{\partial S}{\partial D^{\alpha}q_v}\Bigg) = \pi_{\mu} + H^{\pi}_{\mu} = 0 \, . \quad (15c)$$

Here, the fractional Hamilton's principle function is written as

$$S = S(D^{\alpha-1}q_a, D^{\alpha-1}q_{\mu}, D^{\alpha}q_a, D^{\alpha}q_{\mu}, t) \ , \quad (16)$$

and we can define

$$p_a = \frac{\partial S}{\partial D^{\alpha-1}q_a}, \ p_{\mu} = \frac{\partial S}{\partial D^{\alpha-1}q_{\mu}} \ , \ \pi_a = \frac{\partial S}{\partial D^{\alpha}q_a}, \ \pi_{\mu} = \frac{\partial S}{\partial D^{\alpha}q_{\mu}} \text{ and } p_{\circ} = \frac{\partial S}{\partial t}.$$

Thus, the equations of motion in fractional form can be written as total differential equations as follows [14]:

$$dD^{\alpha-1}q_a = \frac{\partial H'_{\circ}}{\partial p_a}dt + \frac{\partial H'^{p}_{\mu}}{\partial p_a}dD^{\alpha-1}q_{\mu} + \frac{\partial H'^{\pi}_{\mu}}{\partial p_a}dD^{\alpha}q_{\mu}, \quad (17a)$$

$$dD^{\alpha}q_a = \frac{\partial H'_{\circ}}{\partial \pi_a}dt + \frac{\partial H'^{p}_{\mu}}{\partial \pi_a}dD^{\alpha-1}q_{\mu} + \frac{\partial H'^{\pi}_{\mu}}{\partial \pi_a}dD^{\alpha}q_{\mu}, \quad (17b)$$

$$-dp_i = \frac{\partial H'_{\circ}}{\partial D^{\alpha-1}q_i}dt + \frac{\partial H'^{p}_{\mu}}{\partial D^{\alpha-1}q_i}dD^{\alpha-1}q_{\mu} + \frac{\partial H'^{\pi}_{\mu}}{\partial D^{\alpha-1}q_i}dD^{\alpha}q_{\mu}, \quad (17c)$$

$$-d\pi_i = \frac{\partial H'_{\circ}}{\partial D^{\alpha}q_i}dt + \frac{\partial H'^{p}_{\mu}}{\partial D^{\alpha}q_i}dD^{\alpha-1}q_{\mu} + \frac{\partial H'^{\pi}_{\mu}}{\partial D^{\alpha}q_i}dD^{\alpha}q_{\mu}, \quad (17d)$$

If the total derivative of equation (15) is zero [3]

$$dH'_{\circ} = 0; \ \ dH'^{p}_{\mu} = 0; \ \ \ dH'^{\pi}_{\mu} = 0. \quad (18)$$

This means that equations (17) are integrable and the rank of Hessian matrix is $N-R$ . Because of constraints, the degrees of freedom are reduced from $N$ to $N-R$, thus, the constraints reduce the canonical phase space coordinates from $\{D^{\alpha-1}q_i, p_i, D^{\alpha}q_i, \pi_i\}$ to $\{D^{\alpha-1}q_a, p_a, D^{\alpha}q_a, \pi_a\}$.

Eqs. (15) can be written as

$$\frac{\partial S}{\partial t}+H_{\circ}(D^{\alpha-1}q_i,D^{\alpha}q_i,p_a,\pi_a)=0\,; \tag{19a}$$

$$\frac{\partial S}{\partial D^{\alpha-1}q_{\mu}}+H_{\mu}^{p}(D^{\alpha-1}q_i,D^{\alpha}q_i,p_a,\pi_a)=0\,; \tag{19b}$$

$$\frac{\partial S}{\partial D^{\alpha}q_{\mu}}+H_{\mu}^{\pi}(D^{\alpha-1}q_i,D^{\alpha}q_i,p_a,\pi_a)=0\,. \tag{19c}$$

## 3. Hamilton-Jacobi Function and Quantization using WKB approximation

Following refs. [6-8] for investigating the Hamilton-Jacobi function and quantization using the WKB approximation for higher-order singular Lagrangian systems, the fractional Hamilton-Jacobi function can be written as

$$\begin{aligned}S(D^{\alpha-1}q_a,D^{\alpha-1}q_{\mu},D^{\alpha}q_a,D^{\alpha}q_{\mu},t)= f(t)+W_a(D^{\alpha-1}q_a,E_a)+W_a'(D^{\alpha}q_a,E_a,E_a')+\\ f_{\mu}(D^{\alpha-1}q_{\mu})+f_{\mu}'(D^{\alpha}q_{\mu})+A\end{aligned}\,. \tag{20}$$

In this case, we would like to find a solution can be written in a separable form, thus, we can guess a general solution for Eqs. (19) in the form of Eq. (20).

Here, the resulting equation for $f(t)$ has the solution $f(t)=-\sum_{a=1}^{N-R}E_a't$, $E_a'$ are the (N-R) constants of integration and $A$ is some other constant; $D^{\alpha-1}q_{\mu}$ and $D^{\alpha}q_{\mu}$ are treated as independent variables, just as the time $t$ and the remaining functions $W_a(D^{\alpha-1}q_a,E_a)$, $W_a'(D^{\alpha}q_a,E_a,E_a')$, $f_{\mu}(D^{\alpha-1}q_{\mu})$, and $f_{\mu}'(D^{\alpha}q_{\mu})$ are the time-independent Hamilton-Jacobi function.

Once we have found the Hamilton-Jacobi function $S$, the equations of motion can be obtained using the canonical transformations [15], as follows:

$$\eta_a=\frac{\partial S}{\partial E_a'}\,; \tag{21a}$$

$$\lambda_a=\frac{\partial S}{\partial E_a}\,; \tag{21b}$$

$$p_i=\frac{\partial S}{\partial D^{\alpha-1}q_i}\,; \tag{21c}$$

$$\pi_i = \frac{\partial S}{\partial D^{\alpha} q_i}, \tag{21d}$$

The number of $\eta_a$ and $\lambda_a$ are equal to the rank of the Hessian matrix, N-R and $\eta_a$ and $\lambda_a$ are constants and can be determined from the initial conditions.

Because of constrained, the rank of the Hessian matrix is N-R for singular Lagrangian systems. Thus, the fractional wave function $\Psi$ for these systems can be written as [6- 8].

$$\Psi(D^{\alpha-1}q_a, D^{\alpha-1}q_{\mu}, D^{\alpha}q_a, D^{\alpha}q_{\mu}, t) = \left[\prod_{a=1}^{N-R} \psi_{0a}(D^{\alpha-1}q_a)\phi_{0a}(D^{\alpha}q_a)\right] \times \exp\frac{iS(D^{\alpha-1}q_a, D^{\alpha-1}q_{\mu}, D^{\alpha}q_a, D^{\alpha}q_{\mu}, t)}{\hbar}, \tag{22a}$$

$$\psi_{0a} = \frac{1}{\sqrt{p(D^{\alpha-1}q_a)}}; \tag{22b}$$

$$\phi_{0a} = \frac{1}{\sqrt{\pi(D^{\alpha}q_a)}}. \tag{22c}$$

and this wave function satisfies the following conditions:

$$\hat{H}_0'\Psi = 0; \tag{23a}$$

$$H_{\mu}'^{p}\Psi = 0; \tag{23b}$$

$$H_{\mu}'^{\pi}\Psi = 0. \tag{23c}$$

These conditions are obtained when the generalized coordinates and momenta are written as operators:

$$D^{\alpha-1}q_i \to D^{\alpha-1}q_i; \tag{24a}$$

$$D^{\alpha}q_i \to D^{\alpha}q_i; \tag{24b}$$

$$p_i \to \hat{p}_i = \frac{\hbar}{i}\frac{\partial}{\partial D^{\alpha-1}q_i}; \tag{24c}$$

$$\pi_i \to \hat{\pi}_i = \frac{\hbar}{i}\frac{\partial}{\partial D^{\alpha}q_i}; \tag{24d}$$

$$p_0 \longrightarrow \hat{p}_0 = \frac{\hbar}{i}\frac{\partial}{\partial t}. \tag{24e}$$

.

## 4. Examples

**Example1:** We will discuss a fractional second-order regular Lagrangian has the form

$$L=\frac{1}{2}\left((D^{2\alpha}q)^2-(D^{\alpha}q)^2\right). \tag{25}$$

The corresponding generalized momenta are

$$p=-D^{\alpha}q-D^{3\alpha}q\,; \tag{26a}$$

$$\pi=D^{2\alpha}q\,; \tag{26b}$$

The Hamiltonian $H_0$ is calculated as

$$H_{\circ}=pD^{\alpha}q+\frac{1}{2}\pi^2+\frac{1}{2}(D^{\alpha}q)^2\,. \tag{27}$$

The corresponding set of HJPDE, reads

$$H'_{\circ}=p_{\circ}+H_{\circ}=p_{\circ}+pD^{\alpha}q+\frac{1}{2}\pi^2+\frac{1}{2}(D^{\alpha}q)^2\,. \tag{28}$$

We note in this example that there are no primary constraints [1, 2].

Thus, the equations of motion can be calculated as

$$dD^{\alpha-1}q\;=\frac{\partial H'_{\circ}}{\partial p}dt=D^{\alpha}q\,dt, \tag{29a}$$

$$dD^{\alpha}q\;=\frac{\partial H'_{\circ}}{\partial \pi}dt=\pi dt, \tag{29b}$$

$$-dp\;=\frac{\partial H'_{\circ}}{\partial D^{\alpha-1}q}dt=0, \tag{29c}$$

$$-d\pi\;=\frac{\partial H'_{\circ}}{\partial D^{\alpha}q}dt=(p+D^{\alpha}q)dt, \tag{29d}$$

The HJPDE, Eq. (19a), reads

$$H'_{\circ}=p_{\circ}+H_{\circ}=\frac{\partial S}{\partial t}+D^{\alpha}q\frac{\partial S}{\partial D^{\alpha-1}q}+\frac{1}{2}\left(\frac{\partial S}{\partial D^{\alpha}q}\right)^2+\frac{1}{2}(D^{\alpha}q)^2=0\,. \tag{30}$$

Substituting Eq. (20) into (30), we have

$$\frac{\partial f}{\partial t}+D^{\alpha}q\frac{\partial W}{\partial D^{\alpha-1}q}+\frac{1}{2}\left(\frac{\partial W'}{\partial D^{\alpha}q}\right)^{2}+\frac{1}{2}(D^{\alpha}q)^{2}=0 \quad . \tag{31}$$

Since $H_0$ is time-independent, one can write $f(t)=-E't$. Eq. (31) can then be written as

$$-E'+D^{\alpha}q\frac{\partial W}{\partial D^{\alpha-1}q}+\frac{1}{2}\left(\frac{\partial W'}{\partial D^{\alpha}q}\right)^{2}+\frac{1}{2}(D^{\alpha}q)^{2}=0. \tag{32}$$

From Eq. (32), We note that $W$ depends only on $D^{\alpha-1}q$ and $W'$ depends only on $D^{\alpha}q$. This means that

$$\frac{\partial W}{\partial D^{\alpha-1}q}=E; \tag{33a}$$

so that

$$W=D^{\alpha-1}qE. \tag{33b}$$

Substituting Eq. (33b) into (32), one can obtain

$$-E'+D^{\alpha}qE+\frac{1}{2}\left(\frac{\partial W'}{\partial D^{\alpha}q}\right)^{2}+\frac{1}{2}(D^{\alpha}q)^{2}=0. \tag{34}$$

This equation leads to

$$W'(D^{\alpha}q,E,E')=\int\sqrt{2E'+E^{2}-(D^{\alpha}q+E)^{2}}\,dD^{\alpha}q. \tag{35}$$

Thus, the Hamilton-Jacobi function can be obtained as

$$S(D^{\alpha-1}q,D^{\alpha}q,E,E')=-E't+D^{\alpha-1}qE+\int\sqrt{2E'+E^{2}-(D^{\alpha}q+E)^{2}}\,dD^{\alpha}q+A \tag{36}$$

By using Eqs. (21a, 21b), we can obtain the solutions for the generalized coordinates

$$\eta=\frac{\partial S}{\partial E'}=-t+\int\frac{dD^{\alpha}q}{\sqrt{2E'+E^{2}-(D^{\alpha}q+E)^{2}}}; \tag{37a}$$

$$\lambda=\frac{\partial S}{\partial E}=D^{\alpha-1}q+\int\frac{[E-(D^{\alpha}q+E)}{\sqrt{2E'+E^{2}-(D^{\alpha}q+E)^{2}}}dD^{\alpha}q; \tag{37b}$$

Eqs. (37) can be solved, respectively, to give

$$D^{\alpha}q=\sqrt{2E'+E^{2}}\sin(\eta+t)-E; \tag{38a}$$

$$D^{\alpha-1}q=\lambda-E(\eta+t)-\sqrt{2E'+E^{2}}\cos(\eta+t). \tag{38b}$$

The generalized momenta can be determined by using Eqs. (21c,21d), after substituting the result for $D^{\alpha}q$:

$$p = \frac{\partial S}{\partial D^{\alpha-1}q} = E \; ; \tag{39a}$$

$$\pi = \frac{\partial S}{\partial D^{\alpha}q} = \sqrt{2E' + E^2 - (D^{\alpha}q + E)^2} = \sqrt{2E' + E^2}\cos(\eta + t) \,. \tag{39b}$$

The wave function for this example can be determined as

$$\Psi(D^{\alpha-1}q, D^{\alpha}q, t) = \left[\psi_{01}(D^{\alpha-1}q)\phi_{01}(D^{\alpha}q)\right]\exp\frac{iS(D^{\alpha-1}q, D^{\alpha}q, t)}{\hbar} \,, \tag{40}$$

Where

$$\psi_{01} = \frac{1}{\sqrt{p(D^{\alpha-1}q)}} = [E]^{\frac{-1}{2}} \; ; \tag{41a}$$

$$\phi_{01} = \frac{1}{\sqrt{\pi(D^{\alpha}q)}} = \left[2E' + E^2 - (D^{\alpha}q + E)^2\right]^{\frac{-1}{4}} . \tag{41b}$$

and the Hamilton-Jacobi $S$ is given by Eq. (36).

Now let us apply the HJPDE, Eq. (28), to the wave function $\Psi$, after representing the generalized coordinates and momenta as operators:

$$\hat{H}'_{\circ}\Psi = \left[\frac{\hbar}{i}\frac{\partial}{\partial t} + D^{\alpha}q\frac{\hbar}{i}\frac{\partial}{\partial D^{\alpha-1}q} - \frac{\hbar^2}{2}\frac{\partial^2}{\partial (D^{\alpha}q)^2} + \frac{1}{2}(D^{\alpha}q)^2\right]\Psi \,. \tag{42}$$

After some algebra, we have

$$\frac{\hbar}{i}\frac{\partial}{\partial t}\Psi = -E'\Psi \; ; \tag{43a}$$

$$\frac{\hbar}{i}\frac{\partial}{\partial D^{\alpha-1}q}\Psi = E\Psi \; ; \tag{43b}$$

$$-\frac{\hbar^2}{2}\frac{\partial}{\partial (D^{\alpha}q)^2}\Psi = \left[\begin{matrix} -\frac{5\hbar^2}{8}(D^{\alpha}q + E)^2(2E' + E^2 - (D^{\alpha}q + E)^{2)})^{-2} - \\ -\frac{\hbar^2}{4}(2E' + E^2 - (D^{\alpha}q + E)^2)^{-1} + \frac{1}{2}(2E' + E^2 - (D^{\alpha}q + E)^2) \end{matrix}\right]\Psi \; ; \tag{43c}$$

.

Thus, Eq. (42) becomes

$$\hat{H}'_{\circ}\Psi = \begin{bmatrix} -E' + D^{\alpha}qE - \frac{5\hbar^2}{8}(D^{\alpha}q+E)^2(2E'+E^2-(D^{\alpha}q+E)^{2)})^{-2} - \\ \frac{\hbar^2}{4}(2E'+E^2-(D^{\alpha}q+E)^2)^{-1} + \frac{1}{2}(2E'+E^2-(D^{\alpha}q+E)^2) + \frac{1}{2}(D^{\alpha}q)^2 \end{bmatrix}\Psi . \quad (44)$$

Taking the semiclassical limit $\hbar \to 0$ in Eq. (44), we obtain

$$\hat{H}'_{\circ}\Psi = \left[ -E' + D^{\alpha}qE + \frac{1}{2}(2E'+E^2-(D^{\alpha}q+E)^2) + \frac{1}{2}(D^{\alpha}q)^2 \right]\Psi = 0 . \quad (45)$$

**Example 2:** Let us consider the following mathematical singular Lagrangian with two primary first-class constraints has the form

$$L = \frac{1}{2}\left((D^{2\alpha}q)^2 + (D^{2\alpha}q_2)^2\right) + D^{\alpha}q_3 D^{2\alpha}q_3 + D^{\alpha}q_3 D^{\alpha-1}q_3 + D^{\alpha-1}q_2 D^{\alpha}q_2 . \quad (46)$$

The corresponding generalized momenta, Eqs. (7,8) and (12) are

$$p_1 = -D^{3\alpha}q_1 ; \quad (47a)$$

$$p_2 = D^{\alpha-1}q_2 - D^{3\alpha}q_2 ; \quad (47b)$$

$$p_3 = D^{\alpha-1}q_3 = -H_3^p ; \quad (47c)$$

$$\pi_1 = D^{2\alpha}q_1 ; \quad (47d)$$

$$\pi_2 = D^{2\alpha}q_2 ; \quad (47e)$$

$$\pi_3 = D^{\alpha}q_3 = -H_3^{\pi} . \quad (47f)$$

Here, Equations (47c) and (47f) can be written as

$$H_3'^{p} = p_3 - D^{\alpha-1}q_3 = 0 ; \quad (48a)$$

$$H_3'^{\pi} = \pi_3 - D^{\alpha}q_3 = 0 . \quad (48b)$$

and represent as primary constraints [1, 2].

The Hamiltonian $H_0$ is calculated as

$$H_{\circ} = p_1 D^{\alpha}q_1 + (p_2 - D^{\alpha-1}q_2)D^{\alpha}q_2 + \frac{1}{2}(\pi_1^2 + \pi_2^2). \quad (49)$$

The corresponding set of HJPDEs, Eqs. (15), reads

$$H'_\circ = p_\circ + H_\circ = p_1 D^\alpha q_1 + (p_2 - D^{\alpha-1} q_2) D^\alpha q_2 + \frac{1}{2}(\pi_1^2 + \pi_2^2). \tag{50a}$$

$$H_3'^{p} = p_3 - D^{\alpha-1} q_3 = 0\,; \tag{50b}$$

$$H_3'^{\pi} = \pi_3 - D^\alpha q_3 = 0\,. \tag{50c}$$

Here, the Poisson brackets $\{H_3'^{p}, H'_\circ\} = 0$, $\{H_3'^{\pi}, H'_\circ\} = 0$ and $\{H_3'^{p}, H_3'^{\pi}\} = 0$. Therefore, there are no secondary constraints, and these are first-class constraints [1, 2].

The corresponding set of HJPDEs, Eqs. (19), reads

$$H'_\circ = p_\circ + H_\circ = \frac{\partial S}{\partial t} + D^\alpha q_1 \frac{\partial S}{\partial D^{\alpha-1} q_1} + D^\alpha q_2 \left( \frac{\partial S}{\partial D^{\alpha-1} q_2} - D^{\alpha-1} q_2 \right) + \frac{1}{2}\left( \frac{\partial S}{\partial D^\alpha q_1} \right)^2 \frac{1}{2}\left( \frac{\partial S}{\partial D^\alpha q_2} \right)^2 = 0$$

; (51a)

$$H_3'^{p} = \frac{\partial S}{\partial D^{\alpha-1} q_3} - D^{\alpha-1} q_3 = 0\,; \tag{51b}$$

$$H_3'^{\pi} = \frac{\partial S}{\partial D^\alpha q_3} - D^\alpha q_3 = 0\,. \tag{51c}$$

the Hamilton-Jacobi function $S$ Eq. (20) can be written as

$$S(D^{\alpha-1} q_1, D^{\alpha-1} q_2, D^{\alpha-1} q_3, D^\alpha q_1 D^\alpha q_2, D^\alpha q_3, t) = f(t) + W_1(D^{\alpha-1} q_1, E_1) + W_2(D^{\alpha-1} q_2, E_2) + W_1'(D^\alpha q_1, E_1, E_1') + W_2'(D^\alpha q_2, E_2, E_2') + f_3(D^{\alpha-1} q_3) + f_3'(D^\alpha q_3) + A \tag{52}$$

Since the Hamiltonian $H_0$ is time-independent and the coordinates $D^{\alpha-1} q_3$ and $D^\alpha q_3$ are treated as independent variables, one can write

$$f(t) = -(E_1' + E_2')t\,.$$

Substituting $S$ into Eq. (51a), we have

$$-E_1' + D^\alpha q_1 \frac{\partial W_1}{\partial D^{\alpha-1} q_1} + \frac{1}{2}\left( \frac{\partial W_1'}{\partial D^\alpha q_1} \right)^2 - E_2' + D^\alpha q_2 \left( \frac{\partial W_2}{\partial D^{\alpha-1} q_2} - D^\alpha q_2 \right) + \frac{1}{2}\left( \frac{\partial W_2'}{\partial D^\alpha q_2} \right)^2 = 0\,. \tag{52}$$

We note that $W_1$ depends only on $D^{\alpha-1} q_1$ and $W_2$ depends only on $D^{\alpha-1} q_2$. We can then write

$$\frac{\partial W_1}{\partial D^{\alpha-1} q_1} = E_1;$$

so that

$$W_1 = D^{\alpha-1}q_1E_1. \tag{53a}$$

and

$$\frac{\partial W_2}{\partial D^{\alpha-1}q_2} - D^{\alpha-1}q_2 = E_2\,;$$

so that

$$W_2 = D^{\alpha-1}q_2E_2 + \frac{1}{2}(D^{\alpha-1}q_2)^2 \quad . \tag{53b}$$

Substituting Eqs. (53) into (52), we get

$$-E_1' + D^{\alpha}q_1E_1 + \frac{1}{2}\left(\frac{\partial W_1'}{\partial D^{\alpha}q_1}\right)^2 - E_2' + D^{\alpha}q_2E_2 + \frac{1}{2}\left(\frac{\partial W_2'}{\partial D^{\alpha}q_2}\right)^2 = 0\,. \tag{54}$$

Separation of variables in this equation yields

$$\frac{1}{2}\left(\frac{\partial W_1'}{\partial D^{\alpha}q_1}\right)^2 + D^{\alpha}q_1E_1 - E_1' = 0\,; \tag{55a}$$

$$\frac{1}{2}\left(\frac{\partial W_2'}{\partial D^{\alpha}q_2}\right)^2 + D^{\alpha}q_2E_2 - E_2' = 0\,. \tag{55b}$$

The solution of these two equations can be determined as

$$W_1'(D^{\alpha}q_1, E_1, E_1') = \int\sqrt{2E_1' - 2D^{\alpha}q_1E_1}\,dD^{\alpha}q_1\,; \tag{56a}$$

$$W_2'(D^{\alpha}q_2, E_2, E_2') = \int\sqrt{2E_2' - 2D^{\alpha}q_2E_2}\,dD^{\alpha}q_2\,. \tag{56b}$$

Using Eq. (51b), we find $f_3(D^{\alpha-1}q_3) = \frac{1}{2}(D^{\alpha-1}q_3)^2$; and using Eq. (51c), we find $f_3(D^{\alpha}q_3) = \frac{1}{2}(D^{\alpha}q_3)^2$.

With these results, the Hamilton-Jacobi function becomes

$$S(D^{\alpha-1}q_1, D^{\alpha-1}q_2, D^{\alpha-1}q_3, D^{\alpha}q_1D^{\alpha}q_2, D^{\alpha}q_3, t) = (-E_1' - E_2')t + D^{\alpha-1}q_1E_1 + \frac{1}{2}(D^{\alpha-1}q_2)^2 + \int\sqrt{2E_1' - 2D^{\alpha}q_1E_1}\,dD^{\alpha}q_1 + \int\sqrt{2E_2' - 2D^{\alpha}q_2E_2}\,dD^{\alpha}q_2 + \frac{1}{2}(D^{\alpha-1}q_3)^2 + \frac{1}{2}(D^{\alpha}q_3)^2 + A \quad . \tag{57}$$

The solutions for the generalized coordinates can be obtained from the transformations Eqs (21a, 21b)

$$\eta_1 = \frac{\partial S}{\partial E_1'} = -t + \int \frac{dD^\alpha q_1}{\sqrt{2E_1' - 2D^\alpha q_1 E_1}} \; ; \tag{58a}$$

$$\eta_2 = \frac{\partial S}{\partial E_2'} = -t + \int \frac{dD^\alpha q_2}{\sqrt{2E_2' - 2D^\alpha q_2 E_2}} \; ; \tag{58b}$$

$$\lambda_1 = \frac{\partial S}{\partial E_1} = D^{\alpha-1} q_1 + \int \frac{D^\alpha q_1}{\sqrt{2E_1' - 2D^\alpha q_1 E_1}} dD^\alpha q_1 \quad ; \tag{58c}$$

$$\lambda_2 = \frac{\partial S}{\partial E_2} = D^{\alpha-1} q_2 + \int \frac{D^\alpha q_2}{\sqrt{2E_2' - 2D^\alpha q_2 E_2}} dD^\alpha q_2 . \tag{58d}$$

One can solve these four equations, we get

$$D^\alpha q_1 = \frac{E_1'}{E_1} - \frac{E_1}{2}(\eta_1 + t)^2 ; \tag{59a}$$

$$D^\alpha q_2 = \frac{E_2'}{E_2} - \frac{E_2}{2}(\eta_2 + t)^2 \qquad ; \tag{59b}$$

$$D^{\alpha-1} q_1 = \lambda_1 + \frac{E_1'}{E_1}(\eta_1 + t) - \frac{E_1}{6}(\eta_1 + t)^3 ; \tag{59c}$$

$$D^{\alpha-1} q_2 = \lambda_2 + \frac{E_2'}{E_2}(\eta_2 + t) - \frac{E_2}{6}(\eta_2 + t)^3 . \tag{59d}$$

Using Eqs. (21c, 21d) to find the other half of the equations of motion and can be determined, after substituting the results for $D^\alpha q_1$ and $D^\alpha q_2$:

$$p_1 = \frac{\partial S}{\partial D^{\alpha-1} q_1} = E_1 ; \tag{60a}$$

$$p_2 = \frac{\partial S}{\partial D^{\alpha-1} q_2} = E_2 + D^{\alpha-1} q_2 = E_2 + \lambda_2 + \frac{E_2'}{E_2}(\eta_2 + t) - \frac{E_2}{6}(\eta_2 + t)^3 \quad ; \tag{60b}$$

$$p_3 = \frac{\partial S}{\partial D^{\alpha-1} q_3} = D^{\alpha-1} q_3 ; \tag{60c}$$

$$\pi_1 = \frac{\partial S}{\partial D^\alpha q_1} = \sqrt{2E_1' - 2D^\alpha q_1 E_1} = -E_1(\eta_1 + t) ; \tag{60d}$$

$$\pi_2 = \frac{\partial S}{\partial D^{\alpha} q_2} = \sqrt{2E_2' - 2D^{\alpha} q_2 E_2} = -E_2(\eta_2 + t); \quad (60e)$$

$$\pi_3 = \frac{\partial S}{\partial D^{\alpha} q_3} = D^{\alpha} q_3, \quad (60f)$$

where $D^{\alpha-1}q_3$ and $D^{\alpha}q_3$ are arbitrary parameters.

Also, we can determine the equations of motion using Eqs. (17) as

$$dD^{\alpha-1}q_1 = D^{\alpha}q_1 dt, \quad (61a)$$

$$dD^{\alpha-1}q_2 = D^{\alpha}q_2 dt, \quad (61b)$$

$$dD^{\alpha}q_1 = \pi_1 dt, \quad (61c)$$

$$dD^{\alpha}q_2 = \pi_2 dt, \quad (61d)$$

$$-dp_1 = 0, \quad (61e)$$

$$-dp_2 = -D^{\alpha}q_2 dt, \quad (61f)$$

$$-dp_3 = -dD^{\alpha-1}q_3, \quad (61g)$$

$$-d\pi_1 = p_1 dt, \quad (61h)$$

$$-d\pi_2 = (p_2 - D^{\alpha-1}q_2)dt, \quad (61i)$$

$$-d\pi_3 = -dD^{\alpha}q_3, \quad (61j)$$

Solving these equations (61), one gets the same results of Eqs. (59, 60).

The next step is to quantize our singular system. First the wave function (22a) for this example can be obtained as

$$\Psi(D^{\alpha-1}q_1, D^{\alpha-1}q_2, D^{\alpha-1}q_3, D^{\alpha}q_1, D^{\alpha}q_2, D^{\alpha}q_3, t) = \left[\psi_{01}(D^{\alpha-1}q_1)\psi_{02}(D^{\alpha-1}q_2)\phi_{01}(D^{\alpha}q_1)\phi_{02}(D^{\alpha}q_2)\right]\exp\frac{iS}{\hbar}; \quad (62)$$

Where

$$\psi_{01}(D^{\alpha-1}q_1) = \frac{1}{\sqrt{p_1(D^{\alpha-1}q_1)}} = \left[E_1\right]^{-\frac{1}{2}}; \quad (63a)$$

$$\psi_{02}(D^{\alpha-1}q_2) = \frac{1}{\sqrt{p_2(D^{\alpha-1}q_2)}} = \left[E_2 + D^{\alpha-1}q_2\right]^{-\frac{1}{2}}; \quad (63b)$$

$$\phi_{01}(D^{\alpha}q_1)=\frac{1}{\sqrt{\pi_1(D^{\alpha}q_1)}}=\left[2E_1'-2D^{\alpha}q_1E_1\right]^{-\frac{1}{4}}; \tag{63c}$$

$$\phi_{02}(D^{\alpha}q_2)=\frac{1}{\sqrt{\pi_2(D^{\alpha}q_2)}}=\left[2E_2'-2D^{\alpha}q_2E_2\right]^{-\frac{1}{4}}. \tag{63d}$$

and the Hamilton-Jacobi $S$ is given by Eq. (57).

Now let us apply the HJPDEs, Eqs. (50), to the wave function $\Psi$, after representing the generalized coordinates and momenta as operators:

$$\hat{H}_\circ'\Psi=\left[\frac{\hbar}{i}\frac{\partial}{\partial t}+D^{\alpha}q_1\frac{\hbar}{i}\frac{\partial}{\partial D^{\alpha-1}q_1}+D^{\alpha}q_2(\frac{\hbar}{i}\frac{\partial}{\partial D^{\alpha-1}q_2}-D^{\alpha-1}q_2)-\frac{\hbar^2}{2}\frac{\partial^2}{\partial(D^{\alpha}q_1)^2}-\frac{\hbar^2}{2}\frac{\partial^2}{\partial(D^{\alpha}q_2)^2}\right]\Psi; \tag{64a}$$

$$\hat{H}_3'^{p}\Psi=\left[\frac{\hbar}{i}\frac{\partial}{\partial D^{\alpha-1}q_3}-D^{\alpha-1}q_3\right]\Psi; \tag{64b}$$

$$\hat{H}_3'^{\pi}\Psi=\left[\frac{\hbar}{i}\frac{\partial}{\partial D^{\alpha}q_3}-D^{\alpha}q_3\right]\Psi. \tag{64c}$$

After some algebra, we have

$$\frac{\hbar}{i}\frac{\partial}{\partial t}\Psi=(-E_1'-E_2')\Psi; \tag{65a}$$

$$\frac{\hbar}{i}\frac{\partial}{\partial D^{\alpha-1}q_1}\Psi=E_1\Psi; \tag{65b}$$

$$\frac{\hbar}{i}\frac{\partial}{\partial D^{\alpha-1}q_2}\Psi=\left[(D^{\alpha-1}q_2+E_2)-\frac{\hbar}{2i}(D^{\alpha-1}q_2+E_2)^{-1}\right]\Psi; \tag{65c}$$

$$-\frac{\hbar^2}{2}\frac{\partial}{\partial(D^{\alpha}q_1)^2}\Psi=\left[-\frac{5\hbar^2}{8}E_1^2(2E_1'-2D^{\alpha}q_1E_1)^{-2}+\frac{1}{2}(2E_1'-2D^{\alpha}q_1E_1)\right]\Psi; \tag{65d}$$

$$-\frac{\hbar^2}{2}\frac{\partial}{\partial(D^{\alpha}q_2)^2}\Psi=\left[-\frac{5\hbar^2}{8}E_2^2(2E_2'-2D^{\alpha}q_2E_2)^{-2}+\frac{1}{2}(2E_2'-2D^{\alpha}q_2E_2)\right]\Psi; \tag{65e}$$

$$\frac{\hbar}{i}\frac{\partial}{\partial D^{\alpha-1}q_3}\Psi=D^{\alpha-1}q_3\Psi; \tag{65f}$$

$$\frac{\hbar}{i}\frac{\partial}{\partial D^{\alpha}q_3}\Psi = D^{\alpha}q_3\Psi \,. \tag{65g}$$

Substituting the results of Eqs. (65) in Eqs. (64), we obtain

$$\hat{H}'_{\circ}\Psi = \begin{bmatrix} -E_1' - E_2' + D^{\alpha}q_1E_1 + D^{\alpha}q_2[(E_2 + D^{\alpha-1}q_2) - \frac{\hbar}{2i}(E_2 + D^{\alpha-1}q_2)^{-1} - D^{\alpha-1}q_2] - \\ \frac{5\hbar^2}{8}E_1{}^2(2E_1' - 2D^{\alpha}q_1E_1)^{-2} + \frac{1}{2}(2E_1' - 2D^{\alpha}q_1E_1) - \frac{5\hbar^2}{8}E_2^2(2E_2' - 2D^{\alpha}q_2E_2)^{-2}) + \\ \frac{1}{2}(2E_2' - 2D^{\alpha}q_2E_2) \end{bmatrix}\Psi$$

; (66a)

$$H_3'^{p}\Psi = \left[\frac{\hbar}{i}\frac{\partial}{\partial D^{\alpha-1}q_3} - D^{\alpha-1}q_3\right]\Psi = \left[D^{\alpha-1}q_3 - D^{\alpha-1}q_3\right]\Psi = 0 \;; \tag{66b}$$

$$H_3'^{\pi}\Psi = \left[\frac{\hbar}{i}\frac{\partial}{\partial D^{\alpha}q_3} - D^{\alpha}q_3\right]\Psi = \left[D^{\alpha}q_3 - D^{\alpha}q_3\right]\Psi = 0 \;\;. \tag{66c}$$

Taking the limit $\hbar \to 0$ in Eq. (66a), we get

$$\hat{H}'_{\circ}\Psi = \begin{bmatrix} -E_1' - E_2' + D^{\alpha}q_1E_1 + D^{\alpha}q_2[E_2 + D^{\alpha-1}q_2 - D^{\alpha-1}q_2] + \\ \frac{1}{2}(2E_1' - 2D^{\alpha}q_1E_1) + \frac{1}{2}(2E_2' - 2D^{\alpha}q_2E_2) \end{bmatrix}\Psi = 0 \,. \tag{67}$$

**Conclusion**

In our work, we have extended the theory proposed for first-order fractional singular Lagrangian systems to second-order Lagrangian. The set of fractional Hamilton-Jacobi partial differential equations are solved for these systems under the conditions that these equations are separable and satisfy the integrability conditions.

The fractional HJPDEs are obtained using the canonical method. In this work, the fractional Hamilton-Jacobi function S is determined in configuration space for these systems. Finding the function S enables us to get the solutions of the equations of motion.

This is followed by determining the appropriate fractional wave function for these systems. In singular systems the constraints become conditions on the wave function to be satisfied in the semiclassical limit, in addition to the Schrödinger equation, in other words, we have achieved that the quantum results are found to be in exact agreement with the classical results. Finally, we have examined two mathematical examples.

**References**

[1] Dirac, P. A. M. 1950. Generalized Hamiltonian Dynamics. **Canadian Journal of Mathematical Physics**, 2, 129-148.

[2] Dirac, P. A. M. 1964. **Lectures on Quantum Mechanics**, Belfer Graduate School of Science, Yeshiva University, New York.

[3] Guler, Y. 1992b. Canonical Formulation of Singular Systems. **IL Nuovo Cimento B**, 107 (10), 1143-1149.

[4] Rabei, E. M., Nawafleh, K. I., and Ghassib, H. B. 2002. Quantization of Constrained Systems Using the WKB Approximation. **Physical Review A**, 66, 024101-024106.

[5] Muslih, S. I. 2001. Path Integral Formulation of Constrained Systems with Singular Higher-Order Lagrangians. **Hadronic Journal**, 24, 713-721.

[6] Rabei, E. M., Hasan, E. H., and Ghassib, H. B (2004). Hamilton-Jacobi Treatment of Constrained Systems with Second-Order Lagrangians, **International Journal of Theoretical** Physics, Vol. 43 N0. (4), 1073-1096.

[7] Rabei, E. M., Hasan, E. H., and Ghassib, H. B. Quantization of Second-Order Constrained Lagrangian Systems Using the WKB Approximation, **International Journal of Geometric methods in Modern Physics 2005**, Vol. 2. P 485-504.

[8] Hasan, E. H., Rabei, E. M., and Ghassib, H. B. Quantization of Higher-Order Constrained Lagrangian Systems Using the WKB Approximation. **International Journal of Theoretical Physics 2004**, Vol. 43 N0. 11 p 2285-2298**.**

[9] E.M. Rabei and M. Alhorani, Quantization of fractional singular Lagrangian systems using WKB approximation, **international Journal of Modern Physics A**, **Vol. 33** (2018), no. 36, 1850222-1-1850222-9.

[10] E. H. Hasan, Path Integral Quantization of Singular Lagrangians using Fractional Derivatives**, International Journal of Theoretical Physics 2020**, Vol. **59 pages 1157–1164**

[11] E. H. Hasan and J. H. Asad, Remarks on Fractional Hamilton-Jacobi Formalism with second-order Discrete Lagrangian Systems, **Journal of Advanced Physics,** Vol. 6, No. 3, P430-433, 2017

[12] E. H. Hasan, On Fractional Solution of Euler-Lagrange Equations with Second-Order Linear Lagrangians **". Journal of Advanced Physics,** Vol. 7, No. 1, P110-113, 2018

[13] S.G. Samko, A.A. Kilbas and O.I. Marichev, Fractional Integrals and Derivatives: Theory and Applications, Gordon and Breach, 1993.

[14] Pimentel, R.G. and Teixeira 1996. Hamilton-Jacobi Formulation for Singular Systems with Second-Order Lagrangians. **IL Nuovo Ciemento B**, 111, 841-854.

[15] Goldstein, H. 1980. **Classical Mechanics**, 2nd edition, Addison-Wesley, Reading-Massachusetts.